\documentclass[letterpaper, 10 pt, conference]{ieeeconf}  

\IEEEoverridecommandlockouts                              

\usepackage{ourpackages}
\usepackage{tabstackengine}

\def\nochapnum{\par
  \setcounter{chapter}{0}
  \setcounter{section}{0}
  \def\@chapapp{}
}

\def\chapnum{\par
  \setcounter{chapter}{0}
  \setcounter{section}{0}
  \def\@chapapp{Chapter}
}

\newtheorem{fact}{Fact-to-Remember}

\newcommand{\factbox}[3]
{
  \mbox{ }\\
  \centerline{%
    \framebox[0.9\textwidth]{%
      \ifthenelse
      {\equal{#2}{}}
      {\parbox{0.8\columnwidth}{\begin{fact}\label{#1} \mbox{ } \\{#3}\end{fact}}}
      {\parbox{0.8\columnwidth}{\begin{fact}[#2]\label{#1}\mbox{ }\\ {#3} \end{fact}}}
    }
    \mbox{ }\\
  }
}

\newcommand{\mysymb}[2]{
  \noindent
  \begin{tabular}[t]{p{2cm}@{\hspace{1em}:\hspace{1em}}p{9cm}}
    \makebox[2cm][r]{$#1$} & #2\end{ta}\\
}

\newcommand{\bd}{\begin{displaymath}}
    \newcommand{\ed}{\end{displaymath}}

\newcommand{\weg}[1]{}

\newcommand{\identm}[1]{\boldsymbol{\mathbf{\mathit{I}}}_{#1}}

\newcommand{\htrans}[2]%
{\vphantom{T}_{#2}^{\hspace*{0.1ex}#1}%
{\hspace*{-0.1ex}\boldsymbol{T}}\vphantom{T}}

\newcommand{\dotT}[2]%
{\vphantom{T}_{#2}^{\hspace*{0.1ex}#1}%
{\dot{\hspace*{-0.2ex}\boldsymbol{T}}}\vphantom{T}}

\newcommand{\ddotT}[2]%
{\vphantom{T}_{#2}^{\hspace*{0.1ex}#1}%
{\ddot{\hspace*{-0.2ex}\boldsymbol{T}}}\vphantom{T}}

\newcommand{\infhtrans}[2]%
{\vphantom{T}_{#2}^{\hspace*{0.1ex}#1}%
{\hspace*{-0.3ex}\boldsymbol{T}}_{\hspace*{-0.8ex}\Delta}\vphantom{T}}

\newcommand{\esthtrans}[2]%
{\sideset{_{#2}^{\hspace*{0.1ex}#1}}{}{\vphantom{T}}%
  {\hspace*{-0.2ex}\mathop{\widehat{\boldsymbol{T}}}}\vphantom{T}}

\newcommand{\rot}[2]%
{\vphantom{R}_{#2}^{\hspace*{0.1ex}#1}%
{\hspace*{-0.2ex}\boldsymbol{R}}\vphantom{R}}

\newcommand{\dotR}[2]%
{\vphantom{R}_{#2}^{\hspace*{0.1ex}#1}%
{\dot{\hspace*{-0.2ex}\boldsymbol{R}}}\vphantom{R}}

\newcommand{\strans}[2]%
{\vphantom{S}_{#2}^{\hspace*{0.2ex}#1}%
{\hspace*{-0.3ex}\boldsymbol{S}}\vphantom{S}}

\newcommand{\dotS} [2]%
{\vphantom{S}_{#2}^{\hspace*{0.3ex}#1}%
{\dot{\hspace*{-0.3ex}\boldsymbol{S}}}\vphantom{S}}

\newcommand{\infstrans}[2]%
{\vphantom{S}_{#2}^{\hspace*{0.2ex}#1}%
{\hspace*{-0.3ex}\boldsymbol{S}}_{\hspace*{-0.3ex}\Delta}\vphantom{S}}

\newcommand{\sproj}[2]%
{\vphantom{P}_{#2}^{\hspace*{0.2ex}#1}%
{\hspace*{-0.3ex}\boldsymbol{P}}\vphantom{P}}

\newcommand{\dotP} [2]%
{\vphantom{P}_{#2}^{\hspace*{0.3ex}#1}%
{\dot{\hspace*{-0.3ex}\boldsymbol{P}}}\vphantom{P}}

\newcommand{\ptrans}[2]%
{\vphantom{M}_{#2}^{\hspace*{0.2ex}#1}%
{\hspace*{-0.3ex}\boldsymbol{M}}\vphantom{M}}

\newcommand{\dotM}[2]%
{\dot{\hspace*{-0.4ex}\boldsymbol{M}}\vphantom{M}^{#2}_{\hspace*{-0.3ex}#1}}

\newcommand{\infptrans}%
{\boldsymbol{M}_{\hspace*{-0.5ex}\Delta}}

\newcount\m \newcount\n
\def\hours{\n=\time \divide\n 60
  \m=-\n \multiply\m 60 \advance\m \time
  \twodigits\n\ :\ \twodigits\m}
\def\twodigits#1{\ifnum #1<10 0\fi \number#1}

\newcommand{\pkg}[1]{\textsc{#1}}

\newcommand*{\Scale}[2][4]{\scalebox{#1}{$#2$}}%
\DeclareMathOperator*{\minimize}{minimize}

\DeclareRobustCommand{\cpluspluslogo}{\hbox{C\hspace{-0.5ex}
                       \protect\raisebox{0.5ex}
                       {\protect\scalebox{0.67}{++}}}}

\newcommand{\transpose}{^\prime}

\title{\LARGE \bf
\pkg{fatrop}: A Fast Constrained Optimal Control Problem Solver for Robot Trajectory Optimization and Control
}

\author{Lander Vanroye, Ajay Sathya, Joris De Schutter and Wilm Decr\'e
\thanks{All authors are with the Department of Mechanical Engineering, KU Leuven, and with Flanders Make@KU Leuven, Leuven, Belgium.
{\tt\small \{lander.vanroye, ajay.sathya, joris.deschutter, wilm.decre\}@kuleuven.be}.}%
}
\begin{document}

\maketitle
\thispagestyle{empty}
\pagestyle{empty}
\setcounter{MaxMatrixCols}{30}  

\begin{abstract}
        Trajectory optimization is a powerful tool for robot motion planning and control.
        State-of-the-art general-purpose nonlinear programming solvers are versatile, handle constraints effectively and provide a high numerical robustness, but they are slow because they do not fully exploit the optimal control problem structure at hand.
        Existing structure-exploiting solvers are fast, but they often lack techniques to deal with nonlinearity or rely on penalty methods to enforce (equality or inequality) path constraints.
        This work presents \pkg{fatrop}: a trajectory optimization solver that is fast and benefits from the salient features of general-purpose nonlinear optimization solvers.
        The speed-up is mainly achieved through the integration of a specialized linear solver, based on a Riccati recursion that is generalized to also support stagewise equality constraints.
        To demonstrate the algorithm's potential, it is benchmarked on a set of robot problems that are challenging from a numerical perspective, including problems with a minimum-time objective and no-collision constraints.
        The solver is shown to solve problems for trajectory generation of a quadrotor, a robot manipulator and a truck-trailer problem in a few tens of milliseconds.
        The algorithm's \cpluspluslogo-code implementation accompanies this work as open source software, released under the GNU Lesser General Public License (LGPL).
        This software framework may encourage and enable the robotics community to use trajectory optimization in more challenging applications.
\end{abstract}
\section{Introduction}
Nonlinear optimal control problems (OCP) are used in a wide range of engineering applications. 
In robotics, they provide a powerful tool for optimizing robot trajectories and controlling robotic systems.
Trajectory optimization involves finding a control and state trajectory that is (locally) optimal in some metric while satisfying certain equality and inequality path constraints.
This metric, or objective, can, for example, be related to total execution time, safety, energy use or user comfort, while the stagewise path constraints can enforce state and control limits, as well as encode a robot task, such as a path-following task with bounded allowable deviation for a robot end effector.
Trajectory optimization has been applied to many robotics applications such as drone racing \cite{foehn2021time,bos2022multistage}, airborne wind energy systems \cite{horn2013numerical}, legged locomotion \cite{mastalliquadruped} and collision-free robot motion planning \cite{schulman2014motion}.
Apart from trajectory optimization, nonlinear OCPs are building blocks of model predictive control (MPC), which is a control strategy that has gained widespread popularity in the robotics community because of its ability to control complex, underactuated and highly-dynamical systems.
MPC's potential in robotics has been demonstrated in, for example, the control of autonomous racing cars \cite{liniger2015racing} and legged robots \cite{GrandiaMPClegged}.

The computational efficiency of the nonlinear OCP solver is of critical importance for many applications, particularly in reactive sensor-based control applications, where robots have to rapidly adapt their actions to respond to disturbances or unpredictable changes in environment.
Apart from being fast, the solver should also provide a high robustness due to the numerically challenging nature of the considered optimization problems.

We formulate the constrained optimal control problem (COCP) as follows:
\begin{subequations} \label{eq:COCP} 
        \begin{align}
                \minimize_{\substack{\mathbf{x}_k, \mathbf{u}_k,\mathbf{x}_K}} \quad & l_K(\mathbf{x}_K) + \sum_{k=0}^{K-1} l_k(\mathbf{u}_k, \mathbf{x}_k)          \\
                \text{subject to} \quad
                                                                                     & \mathbf{x}_{k+1} = \mathbf{f}_k(\mathbf{u}_k, \mathbf{x}_k),                  \\
                                                                                     & \mathbf{L}_k \leq \mathbf{g}_k(\mathbf{u}_k, \mathbf{x}_k) \leq \mathbf{U}_k, \\
                                                                                     & \mathbf{L}_K\leq \mathbf{g}_K(\mathbf{x}_{K})  \leq \mathbf{U}_K,             \\ \label{eq:COCPeq1}
                                                                                     & \mathbf{h}_k(\mathbf{u}_k, \mathbf{x}_k) =\mathbf{0},                         \\ \label{eq:COCPeq2}
                                                                                     & \mathbf{h}_K(\mathbf{x}_K) = \mathbf{0},
        \end{align}
\end{subequations}
for $k = 0, 1, \dots, K-1$, with $K$ the horizon length, $\mathbf{x}_k$ the state variables, $\mathbf{u}_k$ the control or input variables, $l_k$, $\mathbf{f}_k$, $\mathbf{g}_k$ and $\mathbf{h}_k$ the functions representing the stage cost, the discrete dynamics, the inequality and equality path constraints, respectively.
Note that, because of the general stagewise equality constraints \eqref{eq:COCPeq1}-\eqref{eq:COCPeq2}, this is a broader formulation than common in MPC where, usually, only one stagewise equality constraint is used, namely a constraint at the first time step that fixes the full state vector.
An important problem class, incorporated by this COCP formulation is the boundary value problem where the initial and terminal states are fixed.
The formulation \eqref{eq:COCP} also directly supports moving horizon estimation problems and
furthermore, using helper states, the formulation can encode OCPs with minimum total time objective, multi-stage problems with unknown stage durations as well as periodic systems.

Apart from a difference in formulation, nonlinear MPC and trajectory optimization problems have different solver requirements.
In MPC, approximate optimal control inputs have to be computed within a specified sampling time, due to hard real-time requirements.
An estimate of the solution is available from the previous control step, and only a local neighborhood of this approximate solution has to be explored, which is exploited in MPC-solvers to achieve higher sampling times \cite{verschueren2022acados}.
In trajectory optimization, in contrast, the optimization problem has to be solved to a high accuracy and an estimate of the solution is often not available.
This requires advanced nonlinear programming techniques that provide robustness to cope with the numerically challenging nature of these problems.
General-purpose nonlinear optimization algorithms, such as \pkg{ipopt} \cite{wachter2006implementation}, \pkg{snopt} \cite{gill2005snopt} and \pkg{knitro} \cite{knitro} provide this numerical robustness but they are slow because they do not fully exploit the COCP problem structure at hand.

\subsection{Related Work}
Numerous algorithms tailored for solving nonlinear optimal control problems have been developed, which can be divided into two families: methods that are based on techniques from numerical optimization, such as Sequential Quadratic Programming (SQP), and methods inspired by numerical optimal control, namely  Differential Dynamic Programming (DDP).
For the former family of algorithms techniques that can be employed to exploit the OCP structure include the Riccati Recursion \cite{frison2020hpipm} and specialized condensing approaches \cite{frison2013fast,andersson2013condensing,frison2016pcond}.
A collection of structure-exploiting SQP-based algorithms are implemented in the state-of-the-art MPC framework \pkg{acados} \cite{verschueren2022acados}.
DDP methods, the second family of algorithms, are inspired by the Linear Quadratic Regulator (LQR), which is a well-known concept in linear control.
The computation of the search direction consists of two passes.
In the backward pass, the nonlinear OCP is approximated around the current iterate by a linear quadratic approximation.
This OCP approximation is then used to build an LQR that computes the locally optimal feedback law for each time step.
In the forward pass, the obtained optimal (state) feedback laws and discrete dynamics are used alternately to compute the next iterate.
This procedure is repeated until the problem is solved to the specified precision.
Commonly, these methods use a sum-of-squares objective and neglect the second-order effects of the dynamics and path constraints in the Lagrangian Hessian function.
This results in a Gauss-Newton DDP variant which is referred to as the iterative Linear Quadratic Regulator (iLQR) \cite{giftthaler2018family, todorov2005generalized}.
State-of-the-art trajectory optimization frameworks, for example \pkg{crocoddyl} \cite{mastalli20crocoddyl} and \pkg{altro} \cite{howell2019altro}, implement modified versions of iLQR called feasibility-driven differential dynamic programming (FDDP) and Augmented Lagrangian iLQR (AL-iLQR), respectively. These modifications permit initializing the solver with dynamically infeasible trajectories similar to the multiple shooting method and improve numerical robustness compared to the vanilla iLQR \cite{todorov2005generalized}.
A concise overview of the salient features of the mentioned frameworks is given in Table \ref{tabel:solver_comparison}.
Numerical optimization algorithms have a well-developed set of strategies to solve nonlinear optimization problems with a high robustness.
A difficulty with DDP-based algorithms is that many of these strategies cannot be straightforwardly applied to these methods, and handling nonlinear path constraints robustly is still a developing area \cite{jallet2022constrained}.
A number of challenges that arise in developing DDP-based methods for trajectory optimization have been addressed in prior work.
A first challenge is the treatment of exact Hessian information.
While iLQR methods avoid the, possibly costly, evaluation of second-order derivatives of the dynamics and path constraints, it is often useful to incorporate exact Hessian information.
This is because Gauss-Newton approximations are limited to sum-of-squares objectives, or convex-over-nonlinear objectives in a generalized Gauss-Newton setting \cite{generalized_gn}.
Additionally, exact Hessian information can drastically lower the number of required iterations and increase the solver's numerical robustness.
A difficulty is that line-search algorithms based on exact Hessian DDP, as opposed to iLQR, require Hessian regularization to ensure some descent properties for the filter line search method.
Hessian regularization in exact Hessian DDP was addressed by adding a scaled identity matrix to the full space Hessian \cite{liao1991convergence,pavlov2021interior}.
A second challenge is the treatment of stagewise equality and inequality constraints.
Many DDP methods rely on penalty methods such as quadratic penalty methods (quad pen), augmented Lagrangian methods (ALM) \cite{howell2019altro} and relaxed barrier methods \cite{hauser2006barrier}.
Penalty methods convert constrained OCPs into unconstrained ones, which can be solved with conventional DDP methods.
Some drawbacks of penalty methods are that they require hand tuning, not all are exact, possibly return a local minimizer that does not satisfy the constraints, can suffer from ill-conditioning or can require many iterations.
Inequality constraints that explicitly depend on control variables were implemented with active-set based methods \cite{murray1979constrained,yakowitz1986stagewise,xie2017differential,lin1991differential}.
Inequality constraints were implemented `directly' using a primal-dual interior-point method \cite{pavlov2021interior}.
For state equality constraints with a relative degree of one, a projection method was proposed \cite{giftthaler2017projection}.
A major shortcoming is the lack of existence of a nonlinear optimization solver, tailored for trajectory optimization, that supports arbitrary path constraints, without using penalty methods.
\begin{table}[h]
        \centering
        \Scale[0.9]{ 
                \begin{tabular}{|l|cccc|}
                        \hline
                                             & \pkg{fatrop}  & \pkg{crocoddyl} & \pkg{altro} & \makecell{\pkg{acados}\&\\\pkg{hpipm}}  \\
                        \hline
                        license              & LGPL                     & BSD3            & MIT         & BSD2          \\
                        implementation lang  & C++                      & C++             & \pkg{julia} & C             \\
                        linear algebra       & \pkg{blasfeo}            & \pkg{eigen}     & \pkg{julia} & \pkg{blasfeo} \\
                        dynamics handling    & direct MS                & FDDP            & AL-iLQR     & direct MS     \\
                        stagewise inequality & yes                      & box             & yes         & yes           \\
                        stagewise equality   & yes                      & yes & yes       & no                          \\
                        inequality handling  & NLP-IP                   & quad pen        & ALM         & SQP\&QP-IP    \\
                        equality handling    & \makecell{Newton-\\Lagrange}  & quad pen        & ALM         &  -\\
                        \hline
                \end{tabular}
        }
        \caption{\small Qualitative comparison of different structure-exploiting OCP solvers \pkg{fatrop}, \pkg{acados}, \pkg{altro} and \pkg{crocoddyl}.}
        \label{tabel:solver_comparison}
\end{table}

\subsection{Contribution}
The major contribution of this work is the development of a constrained nonlinear optimal control solver that fast and achieves a high numerical robustness.
A high performance is achieved by exploiting the COCP structure, mainly by the integration of a specialized linear solver introduced in \cite{generalizationriccati}.
A high numerical robustness is obtained by implementing advanced nonlinear programming techniques.
The algorithm is inspired \pkg{ipopt} and handles equality and inequality path constraints in the same way.
The approach uses the multiple shooting formulation which naturally enables initialization from any, possibly infeasible, solution estimate.
Furthermore, the the potential of the approach is demonstrated on a number of benchmark problems from a varying challenging nature.
Finally, the \pkg{fatrop} open source software package that provides an efficient \cpluspluslogo-implementation of the proposed algorithm is released under the LGPL.
\subsection{Outline}
The remainder of this paper is organized as follows.
Section \ref{sec:notation} introduces the notation used throughout this paper.
and discusses preliminaries on nonlinear primal-dual interior point algorithms and direct shooting formulations.
Section \ref{sec:implementation} describes the implementation of the proposed algorithm and how the optimal control problem structure is exploited.
Section \ref{sec:benchmarks} introduces the considered benchmark problems.
Section \ref{sec:implementationdetails} deals with some details on how these problems are translated into optimal control problems and how the benchmark is performed.
Finally, Section \ref{sec:results} presents the results of the benchmark while Section \ref{sec:conclusion} concludes the paper.
Readers who are less familiar with or less interested in the details of the numerical algorithms can skip Sections I-III and jump immediately to Section IV to appreciate the improved performance of the presented solver compared to existing solvers when applied to trajectory optimization problems in robot motion planning and control.
\section{Notation and Preliminaries} \label{sec:notation}
\subsection{Notation}
Vectors are denoted by bold lower case characters, matrices by bold capital characters and scalars by lower case characters.
The symbol $\identm{}$ represents an identity matrix, $\mathbf{0}$ a vector of zeroes and $\mathbf{e}$ a vector of ones,
A diagonal matrix of a vector $\mathbf{x}$ is denoted by $\mathbf{X} = \text{diag}(\mathbf{x})$.
The transpose of matrix $\mathbf{A}$ is denoted by $\mathbf{A} \transpose$.
The following notation represents a linear system $\mathbf{A} \mathbf{x} = -\mathbf{b}$:
       $$  
        \begin{aligned}
                \begin{bNiceMatrix}[first-row]
                        \mathbf{x}                                                            \\
                        \mathbf{A}  & \mathbf{b} \\
                        \CodeAfter
                        \tikz \draw[densely dotted] (1-|2)--(2-|2) ;
                \end{bNiceMatrix}
        \end{aligned}.
       $$ 
The gradient and Hessian of a scalar function $f$ with respect to $\mathbf{x}$ are denoted by $\nabla _\mathbf{x}f$ and $\nabla^2_{\mathbf{xx}}f$, respectively. 
The Jacobian of a vector function $\mathbf{g}$ is denoted by $\mathbf{J}_\mathbf{g}$.
The indices $i$, $k$, and $j$ refer to a specific inequality, time step and interior point outer iteration, respectively.
Blank submatrices in a matrix are structurally zero.
\subsection{Primal-Dual interior-point algorithm}
In this section a concise introduction to primal-dual interior-point methods is given.
\begin{algorithm}[h]\small \label{alg}
        \renewcommand{\thealgorithm}{}
        \begin{algorithmic}[1]
                \While{no convergence of full problem}
                \While{no convergence of barrier subproblem}
                \State compute PD system at current iterate
                \If{not reduced Hessian positive definite}
                \State modify full space Hessian by adding scaled identity matrix
                \EndIf
                \State find search direction (solve primal-dual system)
                \State find step size (line search) and compute next iterate
                \EndWhile
                \State decrease barrier parameter $\mu_j$
                \EndWhile
        \end{algorithmic}
        \caption{\small sketch of the primal-dual interior-point method}
\end{algorithm}

Without loss of generality we consider a slack variable formulation optimization problem of the form:
\begin{subequations} \label{eq:general_opt} 
        \begin{align}
                \minimize_{\mathbf{x}, \mathbf{s}} & \quad {f}(\mathbf{x})                                  \\
                \text{subject to}                   & \quad \mathbf{h}(\mathbf{x}) = \mathbf{0},             \\
                                              & \quad \mathbf{g}(\mathbf{x}) - \mathbf{s} = \mathbf{0}, \\
                                              & \quad  \mathbf{s} \geq \mathbf{0},                      
        \end{align}
\end{subequations}
where $\mathbf{x}$ and $\mathbf{s}$ are the decision and slack variables, respectively.
The slack inequality constraint is implemented by introducing a barrier term to the objective function:
\begin{subequations} \label{eq:general_barrier} 
        \begin{align}
                \minimize_{\mathbf{x}, \mathbf{s}} & \quad {f}(\mathbf{x}) -\mu_j \sum_i \log(s_i)            \\
                \text{subject to}                   & \quad \mathbf{h}(\mathbf{x}) = \mathbf{0},              \\
                                              & \quad \mathbf{g}(\mathbf{x}) - \mathbf{s} = \mathbf{0}. 
        \end{align}
\end{subequations}
After introduction of the additional variables $z_i = \mu_j/s_i$ and some straightforward algebraic manipulation, the first-order necessary optimality conditions associated with \eqref{eq:general_barrier} are:
\begin{subequations}\label{eq:primal_dual_eqs} 
        \begin{align}  
                \nabla_{\mathbf{x}} \mathcal{L} = \mathbf{0}, \\
                \nabla_{\mathbf{s}} \mathcal{L} = \mathbf{0}, \\
                \mathbf{h}(\mathbf{x}) = \mathbf{0},                      \\
                \mathbf{g}(\mathbf{x}) - \mathbf{s} = \mathbf{0},         \\
                \text{diag}(\mathbf{z}) \mathbf{s} = \mu_j  \mathbf{e},       \label{eq:primal_dual_eqs_centering}
        \end{align}
\end{subequations}
where the Lagrangian $\mathcal{L}$ is defined as:
\begin{equation}  
        \mathcal{L} := f(\mathbf{x}) + \boldsymbol{\lambda}_\mathbf{h}\transpose \mathbf{h}(\mathbf{x}) + \boldsymbol{\lambda}_\mathbf{g}\transpose (\mathbf{g}(\mathbf{x}) - \mathbf{s}) - \mathbf{z}\transpose \mathbf{s}.
\end{equation}
Note that these conditions match the first-order optimality conditions of the original problem \eqref{eq:general_opt} as the barrier parameter $\mu_j$ decreases to zero.
The last equation \eqref{eq:primal_dual_eqs_centering}, the centering equation, is a perturbed version of the original problem's complementarity condition.
This means that the method can be viewed effectively as a homotopy method.
The primal-dual interior-point method proceeds by applying Newton's method to this nonlinear system of equations.
This results in a linear system, referred to as the primal-dual system, of the form:
\begin{equation} 
        \begin{aligned}  
                \begin{bNiceMatrix}[first-row]
                        \Delta \mathbf{x}                           & \Delta \mathbf{s} & \Delta \boldsymbol{\lambda_\mathbf{h}} & \Delta \boldsymbol{\lambda_\mathbf{g}} & \Delta \mathbf{z}                                \\
                        \nabla ^2_{\mathbf{x}\mathbf{x}} \mathcal{L} &                  & \mathbf{J_h}\transpose            & \mathbf{J_g}\transpose            &                  & \nabla_\mathbf{x}\mathcal{L} \\
                                                                    &                  &                          & -\identm{}                         & -\identm{}                & \nabla_\mathbf{s}\mathcal{L} \\
                        \mathbf{J_h}                                 &                  &                          &                          &                  & \mathbf{h}                   \\
                        \mathbf{J_g}                                 & -\identm{}                 &                          &                          &                  & \mathbf{g} -\mathbf{s}       \\
                                                                    & \mathbf{Z}                 &                          &                          & \mathbf{S}                & \mathbf{S} \mathbf{z} - \mu_j \mathbf{e}
                        \CodeAfter
                        \tikz \draw[densely dotted] (1-|6)--(6-|6) ;
                \end{bNiceMatrix}
        \end{aligned}.
\end{equation}
Elimination of $\Delta \mathbf{s}$, $\Delta \boldsymbol{\lambda_g}$ and $\Delta \mathbf{z}$ results in a symmetric indefinite linear system of the form:
\begin{equation} \label{eq:reduced_system}
        \begin{aligned}
                \begin{bNiceMatrix}[first-row]
                        \Delta \mathbf{x}                                                                  & \Delta \boldsymbol{\lambda_h}                                                                      \\
                        \nabla ^2_{\mathbf{x}\mathbf{x}}\mathcal{L} + \mathbf{J_g} \transpose \mathbf{S}^{-1} \mathbf{Z} \mathbf{J_g} & \mathbf{J_h} \transpose            & \boldsymbol{\gamma} \\
                        \mathbf{J_h}                                                                       &                          & \mathbf{h}
                        \CodeAfter
                        \tikz \draw[densely dotted] (1-|3)--(3-|3) ;
                \end{bNiceMatrix}
        \end{aligned},
\end{equation}
with $\boldsymbol{\gamma} =  \nabla_\mathbf{x}\mathcal{L} + \mathbf{J}_\mathbf{g} \transpose(-\lambda_\mathbf{g} -\mathbf{S}^{-1}\left(\mu_j \mathbf{e}  - \mathbf{Z}(\mathbf{g}-\mathbf{s})\right))$.
We will refer to this system as the \textit{reduced primal-dual system} throughout the remainder of this paper.
The iteration steps $\Delta \mathbf{s}$, $\Delta \boldsymbol{\lambda_g}$ and $\Delta \mathbf{z}$ can then be retrieved by the following expressions:
\begin{subequations}
        \begin{align}
                 & \Delta \mathbf{s} = \mathbf{J_g} \Delta \mathbf{x} + \mathbf{g} - \mathbf{s},                        \\
                 & \Delta \boldsymbol{\lambda_g} = - \boldsymbol{\lambda_g} -\mathbf{S}^{-1}(\mu_j \mathbf{e} - \mathbf{Z}\Delta \mathbf{s}), \\
                 & \Delta \mathbf{z} = -\mathbf{z} + \mathbf{S}^{-1} (\mu_j \mathbf{e} - \mathbf{Z} \Delta \mathbf{s}).             
        \end{align}
\end{subequations}
The maximum primal and dual step size is chosen in such a way that the fraction-to-boundary rule is satisfied for every slack and dual bound variable, i.e. $s_{i} + \alpha_{\text{primal}}^{\text{max}} \Delta s_i \geq (1-\mu_j) s_{i}$ and $z_{i} + \alpha_{\text{dual}}^{\text{max}} \Delta z_i \geq (1-\mu_j) z_{i}$.
For a more detailed introduction and analysis of primal-dual interior-point methods we refer to the textbook of \textit{Nocedal \& Wright} \cite{nocedal2006numerical}.
\subsection{Direct Multiple Shooting Formulation}
In direct single shooting, all but the first state variable are eliminated by substituting the (discrete-time) dynamics equations.
Multiple shooting, on the other hand, retains all state variables as decision variables, maintaining the dynamics equations as constraints of the nonlinear program.
The latter approach is known to have superior convergence properties over the former in Newton-type optimization algorithms \cite{albersmeyer2010lifted,giftthaler2018family}.
Additionally, multiple shooting allows for initialization from a dynamically infeasible guess, unlike single shooting.
The block-sparse structure of the primal-dual system, arising from the multiple shooting formulation, is exploited in the algorithm's linear solver.
\section{Implementation} \label{sec:implementation}
The nonlinear programming algorithm is heavily inspired by the primal-dual interior-point algorithm \pkg{ipopt} \cite{wachter2006implementation}, applied to the multiple shooting problem formulation.
In this section an overview of the main implementation features of the proposed algorithm is given.
\subsection*{Filter Line Search Globalization}
A filter line search procedure is used to promote global convergence \cite{fletcher2002nonlinear}.
The advantage of using a filter over a merit function acceptance criterion is that the filter's performance is not heavily dependent on the choice of algorithm parameters, such as the merit function's constraint violation penalty parameter.
\subsection*{Exact Hessian Information and Hessian Regularization}
Exact Hessian information represents a problem better locally than Hessian approximations, hence it can drastically lower the number of required iterations and improve the numerical robustness.
A difficulty with using exact Hessian information over Hessian approximation methods like Gauss-Newton and BFGS is that the exact Hessian is not guaranteed to be positive definite.
It is common in line-search methods to require the reduced Hessian approximation, used for the computation of the search direction, to be positive definite.
This guarantees that the computed search direction satisfies some descent properties for the filter line search criterion.
Note that \pkg{fatrop} requires the \textit{reduced} Hessian, this is the full space Hessian projected on the null-space of the constraint Jacobian, to be positive definite.
This is a weaker requirement than positive definiteness of the full space Hessian.
If the reduced Hessian is not positive definite, the full space Hessian is regularized by adding a multiple of the identity matrix.
\subsection*{Structure-Exploiting Linear Solver}
At each iteration, the search direction is computed by solving the reduced primal-dual system \eqref{eq:reduced_system}.
Usually this is the most time consuming step of the algorithm.
The stagewise structure of the COCP \eqref{eq:COCP} results in a block-sparse structure in the reduced primal-dual system.
This block-sparse structure is the same as the KKT system of an equality-constrained OCP and can be exploited by a Riccati recursion that is generalized to also support stagewise equality constraints \cite{generalizationriccati}.
The reduced primal-dual system structure for a horizon length of two $(K=2)$ is shown below:
\begin{equation*}
        \label{eq:blockKKT}
        \Scale[0.8]{
                \begin{bNiceMatrix}[first-row]
                        \mathbf{x}_2        & \mathbf{v}_2           & \boldsymbol{\pi}_2         & \mathbf{u}_1     & \mathbf{x}_1           & \boldsymbol{\lambda}_1         & \boldsymbol{\pi}_1         & \mathbf{u}_0     & \mathbf{x}_0           & \boldsymbol{\lambda}_0         &     \\
                        \mathbf{Q}_2        & \mathbf{H}_2\transpose & -\identm{}    &         &               &                   &               &         &               &                   & \mathbf{q}_2 \\
                        \mathbf{H}_2        &               &               &         &               &                   &               &         &               &                   & \mathbf{h}_2 \\
                        -\identm{} &               &               & \mathbf{B}_1     & \mathbf{A}_1           &                   &               &         &               &                   & \mathbf{b}_1 \\
                                   &               & \mathbf{B}_1\transpose & \mathbf{R}_1     & \mathbf{S}_1\transpose & \mathbf{H}_{1,u}\transpose &               &         &               &                   & \mathbf{r}_1 \\
                                   &               & \mathbf{A}_1\transpose & \mathbf{S}_1     & \mathbf{Q}_1           & \mathbf{H}_{1,x}\transpose & -\identm{}    &         &               &                   & \mathbf{q}_1 \\
                                   &               &               & \mathbf{H}_{1,u} & \mathbf{H}_{1,x}       &                   &               &         &               &                   & \mathbf{h}_1 \\
                                   &               &               &         & -\identm{}    &                   &               & \mathbf{B}_0     & \mathbf{A}_0           &                   & \mathbf{b}_0 \\
                                   &               &               &         &               &                   & \mathbf{B}_0\transpose & \mathbf{R}_0     & \mathbf{S}_0\transpose & \mathbf{H}_{0,u}\transpose & \mathbf{r}_0 \\
                                   &               &               &         &               &                   & \mathbf{A}_0\transpose & \mathbf{S}_0     & \mathbf{Q}_0           & \mathbf{H}_{0,x}\transpose & \mathbf{q}_0 \\
                                   &               &               &         &               &                   &               & \mathbf{H}_{0,u} & \mathbf{H}_{0,x}       &                   & \mathbf{h}_0 \\
                        \CodeAfter
                        \tikz \draw[densely dotted] (1-|11)--(11-|11) ;
                \end{bNiceMatrix}
        },
\end{equation*}
where $\mathbf{x}_k$ and $\mathbf{u}_k$ represent state and input variables, respectively, $\boldsymbol{\pi}_{k}$ represents the dual variables of the discretized dynamics equality constraints, while $\mathbf{v}_K$ and $\boldsymbol{\lambda}_k$ are the dual variables of the stagewise equality constraints.
The block subvectors and matrices appearing in this equation represent the following quantities:
\begin{subequations} \label{eq:blockwise_quantities} \small
        \begin{align}
                 & \begin{bmatrix}
                           \mathbf{H}_{K}
                           \\
                           \mathbf{G}_K
                   \end{bmatrix}
                = \frac{\partial }{\partial \mathbf{x}_K} \begin{bmatrix}
                                                                  \mathbf{h}_K \\ \mathbf{g}_K
                                                          \end{bmatrix},                                                                                                                                  \\
                 &
                \mathbf{Q}_K
                = \frac{\partial ^2  \mathcal{L}}{\partial \mathbf{x}_K \partial \mathbf{x}_K}  + \mathbf{G}_K \mathbf{S}_K^{-1} \mathbf{Z}_K \mathbf{G}_K\transpose,                                                   \\ & \mathbf{q}_K = \nabla \mathbf{x}_K \mathcal{L} + \mathbf{G}_K\transpose(-\boldsymbol{\lambda}_{\mathbf{g}_K} - \mathbf{S}_K^{-1}(\mu_j \mathbf{e} -\mathbf{Z}_K(\mathbf{g}_K-\mathbf{s}_K))),  \\
                 & \begin{bmatrix}
                           \begin{array}{cc}
                                \mathbf{B}_k     & \mathbf{A}_k     \\
                                \mathbf{H}_{k,u} & \mathbf{H}_{k,x}
                        \end{array}
                           \\
                           \mathbf{G}_k
                   \end{bmatrix}
                = \frac{\partial }{\partial \mathbf{w}_k} \begin{bmatrix}
                                                                  \mathbf{f}_k \\ \mathbf{h}_k \\ \mathbf{g}_k
                                                          \end{bmatrix},                                                                                                                   \\
                 & \begin{bmatrix}
                           \mathbf{R}_k & \mathbf{S}_k \transpose \\
                           \mathbf{S}_k & \mathbf{Q}_k            \\
                   \end{bmatrix}
                = \frac{\partial ^2  \mathcal{L}}{\partial \mathbf{w}_k \partial \mathbf{w}_k}  + \mathbf{G}_k \mathbf{S}_k^{-1} \mathbf{Z}_k \mathbf{G}_k \transpose,                                                  \\
                 & \begin{bmatrix} \mathbf{r}_k   \\ \mathbf{q}_k  \end{bmatrix} =  \nabla_{\mathbf{w}_k} \mathcal{L} + \mathbf{G}_k \transpose (-\boldsymbol{\lambda}_{\mathbf{g}_k} - \mathbf{S}_k^{-1}(\mu_j \mathbf{e} - \mathbf{Z}_k(\mathbf{g}_k-\mathbf{s}_k))), 
        \end{align}
\end{subequations}
where $\mathbf{w}_k$ is the concatenation of $\mathbf{u}_k$ and $\mathbf{x}_k$ while $\mathbf{S}_k$ and $\mathbf{Z}_k$ represent the diagonal matrices of the slack variables and dual bound multipliers related to time step $k$.
The used solution scheme only requires a full-rank constraint Jacobian and positive definite reduced Hessian.
These conditions are required by the filter line search anyway.
Moreover, the recursion is used as well to test the positive definiteness of the reduced Hessian at no extra cost, see Step 4 of the Algorithm Sketch.
The computational complexity is linear in the horizon length.
To achieve a high numerical accuracy, iterative refinement is deployed.
For a more detailed overview we refer the reader to \cite{generalizationriccati}.
The \pkg{blasfeo} \cite{frison2018blasfeo} library is used for the linear algebra operations in the algorithm.
This library is optimized for the small-scale matrices that fit in cache memory, appearing in the proposed algorithm.
The kernel routines of this library are optimized for many relevant target CPU architectures and use available CPU capabilities such as \texttt{SIMD} and \texttt{FMA} instruction set extensions.
\subsection*{Stagewise Function Evaluation}
For the problems considered in this paper, the stagewise quantities of the primal-dual system are often the same functions for different time steps.
This means that the code to evaluate these quantities can be re-used.
This way the algorithm has a smaller instruction memory footprint, resulting in a better instruction locality.
Because the code is kept small, it is feasible to apply aggressive compiler optimization levels.
Furthermore, since all block submatrices are independent, they can be evaluated in parallel.
This can be beneficial for problems with expensive function evaluation, for example multi-body problems with forward dynamics.
A preliminary implementation using \pkg{OpenMP} is available in \pkg{fatrop}.
\subsection*{Other Algorithm Features}
Other nonlinear programming features coming from \pkg{ipopt} that are currently implemented in the algorithm are: (1) initialization procedure for the dual and the slack variables, (2) second-order corrections, (3) handling of degenerate constraint Jacobian, (4) watchdog procedure, (5) filter reset heuristic, (6) handling of lower and upper bounds, (7) handling of problems without a strict relative interior and (8) handling of very small search directions.
We refer to the \pkg{ipopt} implementation paper \cite{wachter2006implementation} for a detailed description of these features.
We emphasize that at the moment of writing, the algorithm does not support all features implemented \pkg{ipopt}, such as the feasibility restoration phase, automatic problem scaling and other features that are not described in the original implementation paper \cite{wachter2006implementation}, such as the adaptive barrier parameter update strategy.
Ipopt did not invoke these algorithm features in the experiments of this paper.
\section{Benchmark problems} \label{sec:benchmarks}
Several dynamical systems are considered for which we set up different optimal control problems.
The benchmark problems are of varying dimensions and complexity.
Every problem of the benchmark can be classified either as a model predictive control or as a minimum-time problem.

\noindent \textbf{Cart pendulum: }
an unactuated pendulum is mounted on a cart, which is controlled by a horizontal force, with bounds on the
cart position, velocity and control force. \\ 
\noindent \textit{Model Predictive Control Problem: }
a disturbance is applied while the pendulum is in the upward equilibrium position.
A quadratic objective encodes the task of stabilizing the pendulum while minimizing the total force input over the control horizon. \\ 
\noindent \textit{Swing Minimum Total Time Problem: }
the pendulum starts in the downward configuration at a given cart position.
The goal is to swing the pendulum to an upward position in minimum total time.
Apart from the goal angle, also the translational velocity of the cart and the angular velocity of the pendulum are constrained to be zero at the beginning and end. \\
\textbf{Hanging Chain: }
this nonlinear model predictive control benchmark problem was introduced in \cite{wirsching2006fast}.
The dynamical system consists of a hanging chain of six masses that are connected by springs.
The leftmost mass is rigidly attached to the world and the velocity of the rightmost mass is controlled.
The objective enforces stabilization of the system and the control inputs are limited. We consider both a 2D and a 3D version of this problem.\\
\textbf{Quadrotor: }
\begin{figure}
        \includegraphics[width=0.80\linewidth]{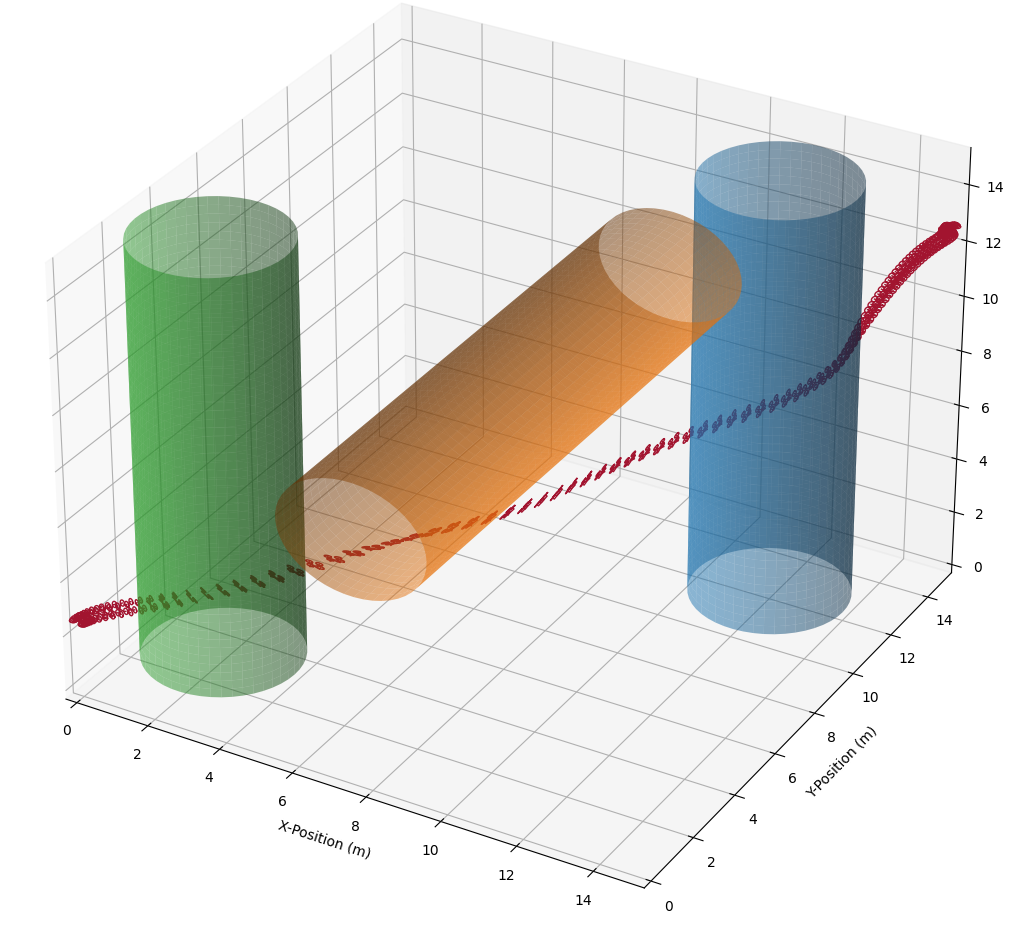} \centering
        \caption{\small Quadrotor no-collision minimum total time problem with three obstacles.}
        \label{fig:drone}
\end{figure}
the orientation is represented by Euler angles. The control inputs are the acceleration in the upward direction of the drone and the Euler angle rates.
The control inputs are limited. \\
\textit{Model Predictive Control Problem: }
the quadrotor starts in horizontal equilibrium position. It is disturbed by a velocity and change in orientation.
The quadratic objective encodes the task to stabilize the quadrotor and go back to the reference position and orientation.
\textit{Point-to-point Minimum Total Time Problem: }
the task is to move the quadrotor from a given initial position and orientation to a given final position and orientation in minimum total time.
The quadrotor has to start and end in equilibrium position. \\ 
\textit{No-collision Minimum Total Time Problem: }
the task is the same as the point-to-point minimum total time problem, but now the quadrotor has to avoid cylindrical obstacles.
There is a variant with a single obstacle and a variant with three obstacles. The latter task is shown in Figure \ref{fig:drone}. \\
\textbf{Seven Degree of Freedom Robot Manipulator: }
the control inputs are the (seven) joint velocity setpoints of the robot.
The robot's joint position and joint velocity limits are taken into account.
\begin{figure}
        \includegraphics[width=0.75\linewidth]{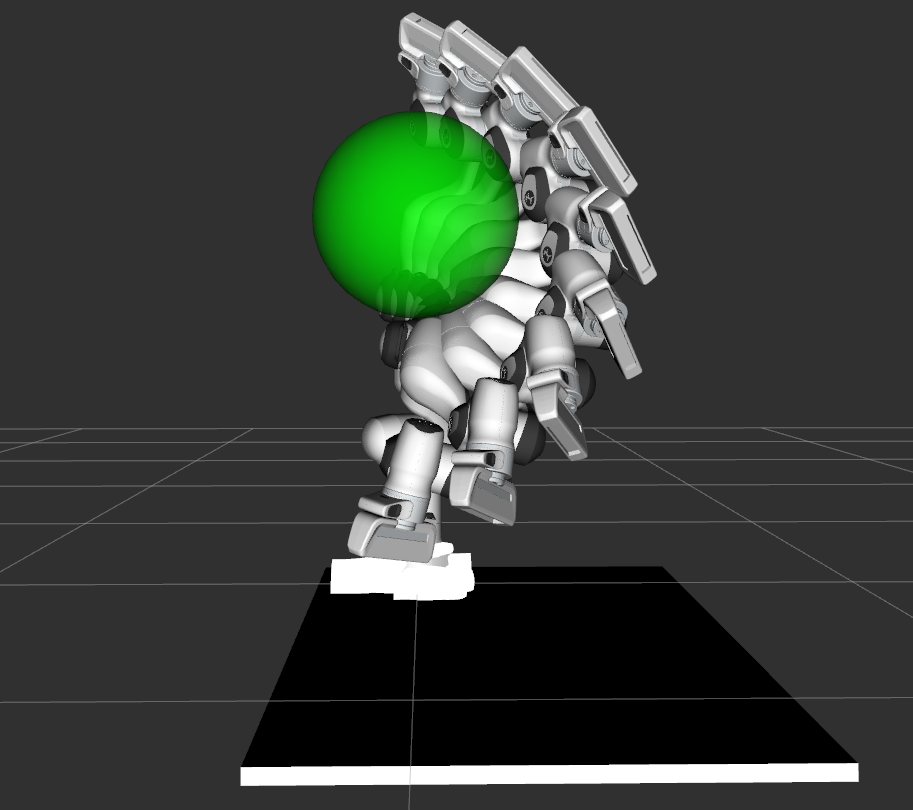} \centering
        \caption{\small Seven degree of freedom robot manipulator with a spherical obstacle.}
        \label{fig:robot}
\end{figure}
The robot has to move from a given position in joint space to a target XYZ position of the end effector in minimum-time.
The robot should avoid collision with a spherical object as shown in Figure \ref{fig:robot}.
The collision model is based on a capsule-based collision model provided by the robot manufacturer. \\
\textbf{Truck with Two Trailers: }
\begin{figure}
\includegraphics[width=0.75\linewidth]{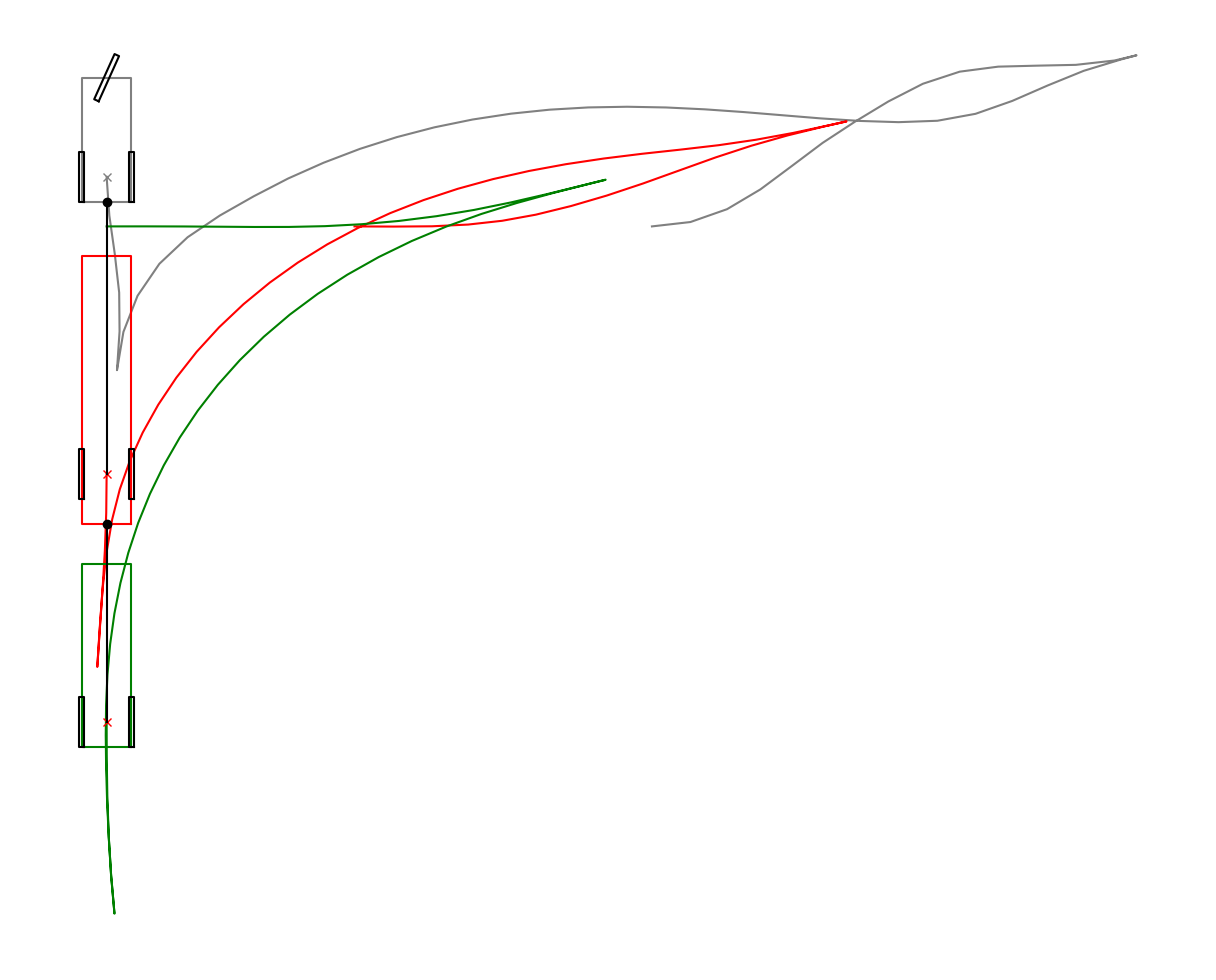} \centering
        \caption{\small Optimal trajectory for truck with two trailers task, the truck and trailer are drawn in the terminal position. The colored lines indicate the trajectory of the coupling of truck and trailers.}
        \label{fig:truck_trailer}
\end{figure}
as indicated in Figure \ref{fig:truck_trailer}, the task involves a truck with two trailers that starts with truck and trailer horizontally aligned and has to park truck and trailer aligned vertically at a given target location in minimum-time.
A kinematic model is used with the control inputs being the steering angular velocities and the velocity of front wheels.


\section{Implementation Details}\label{sec:implementationdetails}
\begin{table}[h]
        \centering
        \Scale[0.9]{ 
                \begin{tabular}{|l|cccccc|}
                        \hline
                        problem                  & $K$ & $n_x$ & $n_u$ & $n_i$ & $n_e[0]$ & $n_e[K]$ \\
                        \hline
                        cart-pendulum mpc        & 25  & 4     & 1     & 0     & 4        & 0        \\
                        cart-pendulum swing*     & 100 & 5     & 1     & 3     & 4        & 3        \\
                        hanging chain 2D mpc     & 25  & 26    & 2     & 2     & 26       & 0        \\
                        hanging chain 2D mpc     & 25  & 39    & 3     & 3     & 39       & 0        \\
                        quadrotor mpc            & 25  & 10    & 4     & 6     & 10       & 0        \\
                        quadrotor p2p*           & 25  & 11    & 4     & 6     & 10       & 8        \\
                        quadrotor one obs*       & 25  & 11    & 5     & 8     & 10       & 8        \\
                        quadrotor three obs*     & 100 & 11    & 7     & 12    & 10       & 8        \\
                        robot manipulator obs*   & 50  & 8     & 18    & 37    & 7        & 3        \\
                        truck with two trailers* & 50  & 8     & 2     & 6     & 5        & 5        \\
                        \hline
                \end{tabular}
        }
        \caption{\small Overview of the problem dimensions. The problems with an asterisk (*) are minimum-time problems. Furthermore, $K$ is the control horizon, $n_x$ the number of states, $n_u$ the number of controls, $n_i$ the number of inequality constraints, and $n_e[0]$ and $n_e[K]$ the number of initial and terminal equality constraints, respectively. 
        The slack variables, necessary for the smooth L1 no-collision constraints formulation, are included in the number of control variables (see Section \ref{sec:implementationdetails}).}
        \label{tab:probdims}
\end{table}
In Table \ref{tab:probdims} the dimensions of each problem are given. \\
Below, we discuss how no-collision constraints and minimum-time problems were handled. \\
\textbf{No-collision constraints. }
in contrast to all other inequality constraints considered in this benchmark set, the no-collision constraints of the Quadrotor and Robot Manipulator were not implemented directly.
This was because we observed in our experiments that the no-collision constraints had difficulties with the combination of the interior points method's strictly feasibility requirements of the slack variables and the nonconvex nature of the considered inequality constraints.
Sometimes the feasibility restoration phase of \pkg{ipopt} was able to overcome this issue but, it sometimes required many iterations, converged to spurious local minima or even failed.
At the moment of writing, \pkg{fatrop} does not implement the feasibility restoration phase, which makes the solver often fail if no-collision constraints are implemented as hard constraints.
As a remedy, we formulated the no-collision constraints as L1-penalized soft constraints (for all solvers), using a smooth reformulation similar to \pkg{TrajOpt}, described in the implementation paper \cite{schulman2014motion}.
The slack variable required by this smooth reformulation is implemented as an auxiliary control variable.
While more efficient handling of slack variables is possible, see for example \pkg{acados} \cite{verschueren2022acados}, this approach was chosen for simplicity of implementation.
We observed that, in our experiments with the \pkg{ipopt} solver, this formulation was more stable and faster than the direct hard constraints implementation of the no-collision constraints.
Because of exactness of the L1-penalty method, the no-collision constraints were always satisfied to specified precision when a large enough penalty parameter was chosen.\\
\textbf{Minimum Total Time Problems. }
To formulate the minimum-time problems in the COCP \eqref{eq:COCP} form, we added an auxiliary state $T_k$, representing the total time variable.
This state was constant over the whole control horizon $(T_{k+1} = T_k)$ and positiveness was enforced by adding the constraint $T_0 \geq 0$.
This $T_k$ variable is then used to determine the integration step $\Delta t  = T_k / K$ of the (uniform) integrator grid. \\
\textbf{Rockit. }
We used \pkg{rockit} \cite{gillis2020effortless} to formulate the OCPs, this is an optimal control problem framework, built on top of \pkg{CasADi} \cite{andersson2019casadi}.
We used the framework's implementation of the Runge-Kutta 4 integrator to transcribe problems' formulations from the continuous time to discrete time.
\pkg{rockit} is interfaced to different solvers, including \pkg{fatrop}, \pkg{ipopt} and \pkg{acados}.
Installation instructions, examples and the code for reproducing the results of the benchmark of this paper made available with this paper \footnote{\url{https://gitlab.kuleuven.be/robotgenskill/fatrop/fatrop_benchmarks}}. \\
\textbf{Benchmarking configuration. } \label{sec:setup}
\pkg{fatrop} was compared to \pkg{ipopt} and, if the problem formulation allowed it, to the \pkg{acados} SQP algorithm.
\pkg{ipopt} is a state-of-the-art general-purpose nonlinear optimization solver while \pkg{acados} is an OCP framework that implements a variety of algorithms with a strong focus on computation speed.
It has some similarities to \pkg{fatrop} as both use \pkg{blasfeo} for linear algebra operations, are able to incorporate exact Hessian information and implement algorithms based on the direct multiple shooting formulation.
The framework's SQP algorithm implements a line search globalization technique.
The problem structure supported by the \pkg{acados} framework is less  general (at the time of this writing) than \pkg{fatrop} (and \pkg{ipopt}) as equality path constraints are not supported, except for constraints that fix the initial state. 
All quantities needed by the algorithms were evaluated from compiled C-code, that was generated using \pkg{CasADi} \cite{andersson2019casadi} SX functions.
We used \pkg{GCC} 9.4.0 with compiler flags \texttt{-Ofast -march=native} for this purpose.
Parallel function evaluation was not implemented for \pkg{ipopt}, so for fairness of comparison, this feature was turned off for both \pkg{acados} and \pkg{fatrop}.
\pkg{ipopt} and \pkg{fatrop} used the same stopping criterion with a tolerance parameter \texttt{tol} of \num{1e-8}. 
We used the default stopping criterion for \pkg{acados}.
\pkg{ipopt} was configured to use the \pkg{ma57} linear solver (sequential), compiled with \pkg{metis} and \pkg{intel mkl}.
\pkg{Blasfeo}, the linear algebra library used by \pkg{fatrop} and \pkg{acados}, was compiled with \texttt{X64\_INTEL\_HASWELL - HP} target.
For \pkg{fatrop} and \pkg{ipopt} we changed the default values of the algorithm parameters \texttt{mu\_init} to \num{1e2} and \texttt{gamma\_theta} to \num{1e-12}, because it benefited the time optimal problems for both solvers.
For \pkg{acados}, \pkg{hpipm} was used as inner QP solver with \texttt{EXACT\_HESSIAN} and \texttt{CONVEXIFY} presets. 
We tried out different numbers of partial condensing steps and took the results for the best performing setting.
We provided discrete-time \pkg{CasADi} integrator expressions, because it resulted in function evaluation that was roughly two to three times faster than providing the continuous-time dynamics differential equations and using the built-in integrator.
Our test machine was a notebook computer equipped with an Intel\textregistered~ Core\texttrademark~ i7-10850H Processor, running Ubuntu 20.04, with  Intel\textregistered~ Turbo Boost\texttrademark~ disabled.
\section{Numerical Results} \label{sec:results}
\pkg{fatrop} and \pkg{ipopt} were able to solve all minimum-time benchmark problems to the specified accuracy, even though no initial solution guess was provided.
A concern might be that the reduced system can become ill-conditioned when iterates come close to the inequality barriers.
We did not observe that this led to numerical issues in our experiments, although many inequality constraints were active at the solution of the minimum-time problems.
The benchmark results for the minimum-time problems are shown in Table \ref{tab:walltimeipopt}. 
Although the iterations for \pkg{fatrop} and \pkg{ipopt} were similar for all problems, they were not exactly the same. This is because of the inexact arithmetic of the linear solvers leading to small step differences that are accumulated over the iterations. In some cases this led to a (small) difference in number of iterations.
The \pkg{MA57} linear solver used in \pkg{ipopt} was usually more accurate than the linear solver used in \pkg{fatrop}, but this did not lead to a significant difference in the number of iterations or robustness.
We refer to the linear solver's paper \cite{generalizationriccati} for a numerical comparison to some general-purpose sparse linear solvers.
\pkg{fatrop} always outperformed \pkg{ipopt}, solving all problems in a few tens of milliseconds to about a hundred milliseconds.
\begin{table}[h] 
        \centering
        \Scale[0.8]{       
                \begin{tabular}{|l|rrr|rrr|} \hline
                        \multirow{2}{*}{problem name} & \multicolumn{3}{c|}{\pkg{fatrop}} & \multicolumn{3}{c|}{\pkg{ipopt}}                                                                                                                                               \\ 
                                                      & \multicolumn{1}{c}{\#it}          & \multicolumn{1}{c}{$t_{\text{tot}}$} & \multicolumn{1}{c|}{$t_\text{FE}$} & \multicolumn{1}{c}{\#it} & \multicolumn{1}{c}{$t_\text{tot}$} & \multicolumn{1}{c|}{$t_\text{FE}$} \\ \hline 
                        cart pendulum swing           & 88                                &  \textbf{35.37}                      &  8.33                              & 81                       &  125.29                            &   11.63                          \\ 
                        quadrotor p2p                & 63                                &  \textbf{10.44}                      &  1.94                              & 63                       &  79.96                             &   3.58                           \\ 
                        quadrotor one obs            & 86                                &  \textbf{17.02}                      &  3.02                              & 79                       &  100.70                            &   3.94                           \\ 
                        quadrotor three obs          & 170                               &  \textbf{135.11}                     &  26.64                             & 169                      &  706.28                            &   34.74                          \\ 
                        robot manipulator obs         & 38                                &  \textbf{35.51}                      &  11.88                             & 38                       &  136.65                            &   54.23                          \\ 
                        truck with two trailers       & 91                                &  \textbf{30.32}                      &  11.80                             & 83                       &  135.63                            &  16.39                           \\ \hline 
                \end{tabular}
        }
        \caption{\small Number of iterations, total and function evaluation wall time in milliseconds (ms) for the minimum time problems with \pkg{fatrop} and \pkg{ipopt}.}
        \label{tab:walltimeipopt}
\end{table}
Only the problems with an MPC formulation were benchmarked with \pkg{acados}, because this solver did not allow the problem formulation of the minimum time problems due to the presence of terminal (or path) equality constraints.
The results for these problems are shown in Table \ref{tab:walltimeipopt}.
\begin{table}[h] 
        \centering
        \Scale[0.8]{        \setlength{\tabcolsep}{2pt}
                \begin{tabular}{|l|rrr|rrr|rrr|} \hline
                        \multirow{2}{*}{problem name} & \multicolumn{3}{c|}{\pkg{fatrop}} & \multicolumn{3}{c|}{\pkg{acados}}    & \multicolumn{3}{c|}{\pkg{ipopt}}                                                                                                                                                                                                             \\ 
                                                      & \multicolumn{1}{c}{\#it}          & \multicolumn{1}{c}{$t_{\text{tot}}$} & \multicolumn{1}{c|}{$t_\text{FE}$} & \multicolumn{1}{c}{\#it} & \multicolumn{1}{c}{$t_\text{tot}$} & \multicolumn{1}{c|}{$t_\text{FE}$} & \multicolumn{1}{c}{\#it} & \multicolumn{1}{c}{$t_\text{tot}$} & \multicolumn{1}{c|}{$t_\text{FE}$} \\ \hline 
                        cart pendulum mpc             & 6                                 & $\textbf{0.54} $                     & 0.15                               & 2                        & 0.59                               & 0.10                               & 6                        & 3.31                               & 0.23                               \\ 
                        quadrotor mpc                & 12                                & $\textbf{1.71} $                     & 0.41                               & 5                        & 5.31                               & 0.39                               & 12                       & 9.61                               & 0.65                               \\ 
                        hanging chain 2D mpc          & 14                                & 7.61                                 & 4.02                               & 4                        & $\textbf{5.22}$                    & 1.93                               & 14                       & 60.48                              & *                                  \\ 
                        hanging chain 3D mpc          & 14                                & 17.62                                & 10.29                              & 3                        & $\textbf{9.76}$                    & 4.51                               & 14                       & 83.96                              & *                                  \\ \hline 
                \end{tabular}
        }
        \caption{\small Number of iterations, total and function evaluation wall time in milliseconds (ms) for the MPC problems with \pkg{fatrop}, \pkg{acados} and \pkg{ipopt}.
                For the results with an asterisk (*), function evaluation time is excluded from the total time because the CasADi virtual machine was used instead of compiled generated C-code.
                For these problems compilation failed due to insufficient available memory on the test machine.} \label{tabel:walltimempc}
\end{table}
\pkg{fatrop} and \pkg{acados} were always faster than \pkg{ipopt}, solving all problems in a few milliseconds. 
The \pkg{acados} algorithm always converged in fewer iterations than \pkg{fatrop}, which explains why its function evaluation time was always lower than \pkg{fatrop}.
The SQP algorithm of \pkg{acados} solves a QP at every (outer) iteration which is computationally more expensive than solving the primal-dual system in a \pkg{fatrop} iteration.
This explains why \pkg{fatrop} was faster for the cart pendulum and quadrotor MPC problems, despite the higher number of iterations.
\pkg{acados}'s inner QP solver \pkg{hpipm} benefits from partial condensing \cite{frison2016pcond}, as the condensed quantities only have to be computed for every outer iteration.
Partial condensing is most beneficial for problems with a large number of states compared to the number of controls ($n_x \gg n_u$).
This, together with the lower number of function evaluations, explains why \pkg{acados} is faster than \pkg{fatrop} for the hanging chain problems.

\section{Conclusion} \label{sec:conclusion}
\pkg{fatrop} is a novel trajectory optimization framework that is heavily inspired by the \pkg{ipopt} algorithm, but achieves a speed-up by exploiting the optimal control structure at hand.
We demonstrated its versatility, numerical robustness and efficiency by solving a variety of benchmark problems.
The \pkg{fatrop} software is available at \url{https://github.com/meco-group/fatrop} under the LGPL.
Future work includes direct handling of no-collision constraints, by means of a specialized feasibility restoration phase, and optimization over Lie manifolds.
Additionally, at this moment the \pkg{fatrop}-\pkg{rockit} interface is limited to single stage problems. Another part of future work is the ability to transcribe multi-stage \pkg{rockit} problems to \pkg{fatrop}.
\section*{\small Acknowledgements}
\small
This result is part of a project that has received funding from the European Research Council (ERC) under the European Union's Horizon 2020 research and innovation programme (Grant agreement: ROBOTGENSKILL No. 788298) and from the Research Foundation Flanders (FWO) (Grant agreement No. G0D1119N).
The authors would like to thank Bastiaan Vandewal and Mathias Bos for providing the dynamical models used for the truck-trailer and quadrotor benchmark problems, respectively.
\bibliographystyle{IEEEtran}
\bibliography{references}

\end{document}